\documentclass[12pt,twoside]{amsart}\usepackage{amssymb}
\usepackage{enumerate}

\newtheorem{thmspec}{\relax}
\newtheorem{theorem}{Theorem}[section]
\newtheorem{thm}[theorem]{Theorem}
\newtheorem{lem}[theorem]{Lemma}

\newtheorem{prop}[theorem]{Proposition}

\theoremstyle{definition}

\theoremstyle{remark}

\numberwithin{equation}{section} \tolerance = 10000

\def \Bbb{\mathbb}

\def\onto{{\kern3pt\to\kern-8pt\to\kern3pt}}

\def\<{\langle}
\def\>{\rangle}
\def\|{{\ |\ }}

\def\onto{\twoheadrightarrow}

\def\-{\underline}

\def\mes{\operatorname{mes}}

\def\End{\operatorname{End}}
\def\arg{\operatorname{arg}}

\def\N{\Bbb N}

\def\R{\Bbb R}
\def\C{\Bbb C}

\def\X{\Bbb X}

\def\<{\langle}
\def\>{\rangle}

\catcode`\@=11
\def\serieslogo@{\relax}
\def\@setcopyright{\relax}
\catcode`\@=12

\title[Boundary cross theorem in dimension $\boldsymbol{1}$ with singularities]
{Boundary cross theorem in dimension $\boldsymbol{1}$ with
singularities}

\begin{document}

\author{Peter Pflug}
\address{Peter Pflug\\
Carl von  Ossietzky Universit\"{a}t Oldenburg \\
Fachbereich  Mathematik\\
Postfach 2503, D--26111\\
 Oldenburg, Germany}
\email{pflug@mathematik.uni-oldenburg.de}

\author{Vi{\^e}t-Anh  Nguy\^en}
\address{Vi{\^e}t-Anh  Nguy\^en\\
Mathematics Section\\
The Abdus Salam international centre
 for theoretical physics\\
Strada costiera, 11\\
34014 Trieste, Italy} \email{vnguyen0@ictp.trieste.it}



\begin{abstract}
Let $D$ and $G$ be copies of the open unit disc in $\C,$  let
  $A$ (resp. $B$) be a  measurable subset of  $\partial D$ (resp.
  $\partial G$),  let $W$ be the  $2$-fold cross $\big((D\cup A)\times B\big)\cup \big(A\times(B\cup G)\big),$
  and let $M$ be a relatively  closed subset of $W.$ 
  Suppose in addition   that
  $A$ and $B$ are of positive one-dimensional Lebesgue measure and that $M$ is fiberwise polar (resp. fiberwise discrete)
  and that $M\cap (A\times B)=\varnothing.$ 
   We  determine the ``envelope of holomorphy"
  $\widehat{W\setminus M}$ of $W\setminus M$ in the sense that any function  locally bounded on $W\setminus M,$
  measurable on $A\times B,$ and separately
  holomorphic
  on $\big((A\times G) \cup (D\times B)\big)\setminus M$ ``extends" to a function
   holomorphic  on $\widehat{W\setminus M}.$
\end{abstract}
\maketitle

\section{Introduction}
The first boundary  cross
theorem
was discovered by  Malgrange--Zerner in the  pioneer  work  \cite{ze}.
Subsequent  results in this direction are obtained by
Komatsu \cite{ko} and Dru\.{z}kowski \cite{dr}.
Recently, Gonchar \cite{go1,go2} has proved  a  more
general result for the one-dimensional case. It should be noted that  
  Airapetyan and Henkin publish  a  version of
 the edge-of-the-wedge theorem for CR manifolds (see  \cite{ah1}  for a brief version
 and  \cite{ah2} for a complete proof).
   Gonchar's result could be deduced from the latter works.
 In  the articles  \cite{pn1, pn2, pn3, pn}, the authors  generalize Gonchar's result
 to the  one-dimensional case with more optimal hypotheses and  to the higher
dimensional case.

On the other hand,  cross theorems with analytic or pluripolar singularities have been developed by
many mathematicians (see, for example, \cite{jp2,jp3,jp4,jp5}
and the references therein). The question naturally arises whether
there exists a mixture of these two types of cross theorems, namely, a  boundary cross theorem  with singularities.

The purpose of this article is to establish  such a theorem in a
simple  but very useful setting:   one-dimensional case with optimal hypotheses in the spirit of our previous work \cite{pn2}.
This  is  our  first  step towards  a general  cross theorem  with singularities  \cite{pn5} (see also \cite{nv1,nv2}).

\smallskip

\indent{\it{\bf Acknowledgment.}} The  paper  was  started during
the  stay of the second  author at the University of Oldenburg in  2006. He  was  supported  by  a  grant from DFG, Az. PF 227/8-2.
   The paper was written while
he was visiting the  Abdus Salam International
Centre
 for Theoretical Physics
in Trieste. He wishes to express his gratitude to these
organizations.

\section{Background and statement of the main result}

First we introduce some notation and terminology. In this article,
 $E$ always denotes the open unit disc in $\C.$ For $a\in\C$ and $r>0,$
 $\Delta_a(r)$ is the disc centered at $a$ with radius $r.$
 Finally,  the one-dimensional
Lebesgue  measure is  denoted by $\mes.$
\subsection{(Sub)harmonic measure}
Let $\Omega\subset \C$ be an open set.  For any function $u:\ \Omega\longrightarrow \R\cup\{-\infty\},$ let
\begin{equation*}
 \hat{u}(z):=
\begin{cases}
u(z),
  & z\in   \Omega,\\
 \limsup\limits_{\Omega\ni w\to z}u(w), & z \in \partial \Omega.
\end{cases}
\end{equation*}
For a set  $A\subset \overline{\Omega}$ put
\begin{equation*}
h_{A,\Omega}:=\sup\left\lbrace u\ :\  u\in\mathcal{SH}(\Omega),\
u\leq 1\ \text{on}\ \Omega,\
   \hat{u}\leq 0\ \text{on}\ A    \right\rbrace,
\end{equation*}
where $\mathcal{SH}(\Omega)$ denotes the cone   of all functions
subharmonic on $\Omega.$

The {\it subharmonic measure} of $A$ relative to $\Omega$ is
the function  $\omega(\cdot,A,\Omega)\in \mathcal{SH}(\Omega)$
defined by
\begin{equation*}
\omega(z,A,\Omega):=h^{\ast}_{A,\Omega}(z),\qquad z\in\Omega,
\end{equation*}
where  $h^{\ast}$ denotes the upper semicontinuous regularization
of  $h.$

If  $A\subset\partial \Omega,$ then
$\omega(\cdot,A,\Omega)$ is also called
 the {\it  harmonic measure} of $A$ relative to $\Omega.$
 In this case,  $\omega(\cdot,A,\Omega)$
is a harmonic function.

We recall  the  following elementary property which  will be  used  several times  later on.
Let  $(A_k)_{k=1}^{\infty}$  be a sequence of measurable subsets of $\partial E$
and  $A$ a  measurable subset of $\partial E$ such that
$
\mes(A_k)>0,$  $ A_k\subset A_{k+1},$  and $  \mes\big( A\setminus \bigcup_{k=1}^{\infty}A_k\big)=0.$
Then
\begin{equation}\label{eq_elementary}
\omega(\cdot,A_k,E)\searrow\omega(\cdot,A,E)\qquad\text{as}\ k\nearrow\infty.
\end{equation}
\subsection{Angular approach regions and locally regular points}
Let $D\subset \C$  be  a  Jordan  domain. Fix a conformal  mapping $\Phi$ from
$D$  onto $E$  which  extends continuously from $\overline{D}$  onto $\overline{E}.$
For   $\zeta\in\partial D$ and  $0<\alpha <\frac{\pi}{2},$   the {\it Stolz region} or
{\it angular approach region} $\mathcal{A}_{\alpha}(\zeta)$ is given by
\begin{equation*}
\mathcal{A}_{\alpha}(\zeta):=
 \left\lbrace   \Phi^{-1}(t):\        t\in E\  \text{and}\ \left\vert
 \arg\left(\frac{\Phi(\zeta)-t}{\Phi(\zeta)}\right)
 \right\vert<\alpha\right\rbrace   ,
\end{equation*}
where  $\arg:\ \C\longrightarrow (-\pi,\pi]$ is  as usual the argument function.


Let  $A\subset\overline D.$ We say that a
point
 $\zeta\in  \overline{ D} $  is a
   {\it  locally regular point relative to}
  $A$
if
\begin{equation*}
\lim\limits_{D\cap \Delta_{\zeta}(r)\ni z\to \zeta} \omega(z,A\cap
\Delta_{\zeta}(r),D\cap \Delta_{\zeta}(r))=0,\qquad r>0.
\end{equation*}
 Obviously, $\zeta\in \overline{A}.$  The set of all locally regular points relative to $A$ is denoted
by $A^{\ast}.$   $A$ is said to be {\it locally regular} if $A\subset
A^{\ast}.$

If $A\subset \partial D$ is  measurable, then it is classical that $\Phi(A^{\ast})$ contains all
density-points of $\Phi(A),$  hence $\mes\Big(\Phi\big(A\setminus(A\cap
A^{\ast})\big)\Big)=0,$  and   $A\cap A^{\ast}$ is  again locally
regular. Moreover, it follows  from (\ref{eq_elementary}) that
\begin{equation}\label{eq_elementary_new}
\omega(\cdot,A\cap A^{\ast},D)=\omega(\cdot,A,D).
\end{equation}

Recall  from Definition 4.8 in \cite{pn2} the following
definition.
A point $\zeta\in\partial D$ is
said to be an {\it  end-point} of an open subset  $\Omega\subset D$ if, for every
$0<\alpha<\frac{\pi}{2},$ there is an open neighborhood
$U=U_{\alpha}$ of $\zeta$ such that $U\cap
\mathcal{A}_{\alpha}(\zeta)\subset \Omega.$ The set of all
end-points of $\Omega$  is denoted by $\End(\Omega).$

We say that a function $f,$  defined in an open subset $\Omega\subset D,$ admits an {\it angular limit}
  $\lambda\in\C$ at a point  $a\in\End(\Omega) $ if
 \begin{equation*}
\lim\limits_{ \mathcal{A}_{\alpha}(a)\cap \Omega\ni z\to a }
f(z)=\lambda,\qquad 0<\alpha<\frac{\pi}{2}.
\end{equation*}

\subsection{Cross and separate holomorphy}

Let $D,G\subset \C $ be two open sets, let
  $A$ (resp. $B$) be a subset of  $\overline{D}$ (resp.
  $\overline{G}$). We define
a {\it $2$-fold cross} $W,$ its {\it  interior} $W^{\text{o}}$  as
\begin{eqnarray*}
W &:=&\X(A,B; D,G)
:=\big((D\cup A)\times B\big)\cup \big(A\times(B\cup G)\big),\\
W^{\text{o}} &:=&\X^{\text{o}}(A,B; D,G)
:=  (D\times B)\cup (A\times  G).
\end{eqnarray*}
For a $2$-fold cross $W :=\X(A,B; D,G)$
define
\begin{equation*}
\widehat{W}=\widehat{\X}(A,B;D,G)
:=\left\lbrace (z,w)\in D\times G:\ \omega(z,A,D)+\omega(w,B,G)  <1
\right\rbrace.
\end{equation*}
Let $M$ be a subset of $W.$  Then the {\it fibers} $M_a$  and $M^b$ are given by
 \begin{equation*}
 M_a:=\{w\in G:\  (a,w)\in M\}\quad  (a\in A); \qquad
 M^b:=\{z\in D:\  (z,b)\in M\}\quad (b\in B).
 \end{equation*}
 We say  that  $M$ possesses  a certain  property in fibers over  $A$  (resp. over  $B$) if
 all  fibers  $M_a$  with $a\in A$ (resp.  all  fibers  $M^b$  with $b\in B$) possess this property.

Suppose that $M$  is  relatively closed  in fibers over  $A$ and  $B.$
We say that a function $f:W\setminus M\longrightarrow \C$ is {\it separately holomorphic}
{\it on $W^{\text{o}}\setminus M $} and write $f\in\mathcal{O}_s(W^{\text{o}}\setminus M),$   if
 for any $a\in A $ (resp.  $b\in B$)
 the function $f(a,\cdot)|_{G\setminus M_a}$  (resp.  $f(\cdot,b)|_{D\setminus M^b}$)
   is holomorphic.

From  now  on we  assume, in addition, that $D$ and $G$ are Jordan domains, and
$A\subset \partial D,$  $B\subset\partial G.$ Then we define the {\it regular part} $W^{\ast}$ relative to $W$ as
\begin{equation*}
 W^{\ast} :=\X(A^{\ast},B^{\ast}; D,G).
 \end{equation*}

Let  $\Omega$ be an open subset of  $D\times G.$
A point  $(a,b)\in A^{\ast}\times G$  (resp.  $(a,b)\in D\times B^{\ast}$)  is  said  to be an {\it end-point}  of $\Omega$
if, for every
$0<\alpha<\frac{\pi}{2},$ there are an open neighborhood
$U=U_{\alpha}$ of $a$ and an open neighborhood $V=V_{\alpha}$ of $b$  such that
\begin{equation*}
\big(U\cap
\mathcal{A}_{\alpha}(a)\big)\times V\subset \Omega\quad\Big(\text{resp.}\ U\times \big(V\cap
\mathcal{A}_{\alpha}(b)\big)\subset \Omega\ \Big).
\end{equation*}
The set of all end-points of $\Omega$  is denoted by $\End(\Omega).$
We say that a function $f:\  \Omega\longrightarrow \C$ admits an {\it angular limit}
  $\lambda\in\C$ at $(a,b)\in  \End(\Omega) $ if
  under  the previous  notation  one  of the following cases occurs:\\
{\bf Case 1:}  $(a,b)\in\ A^{\ast}\times G$ and the  following limits  exist and are equal  to $\lambda$
 \begin{equation*}
\lim\limits_{\Omega\ni (z,w)\to (a,b),\  z \in  \mathcal{A}_{\alpha}(a)}
f(z,w),\qquad 0<\alpha<\frac{\pi}{2};
\end{equation*}
{\bf Case 2:} $(a,b)\in D\times B^{\ast}$ and the  following limits  exist  and are equal  to $\lambda$
 \begin{equation*}
\lim\limits_{\Omega\ni (z,w)\to (a,b),\ w \in  \mathcal{A}_{\alpha}(b)}
f(z,w),\qquad 0<\alpha<\frac{\pi}{2}.
\end{equation*}


For an open set $\Omega\subset\C^k,$ let $\mathcal{O}(\Omega)$ denote
the space of all holomorphic functions on $\Omega.$
A function $f:\  \mathcal{P}\rightarrow \C,$  where $ \mathcal{P}$ is  a topological space,
is said to be {\it locally bounded,} if for every point $p\in  \mathcal{P}$
there exists a neighborhood $U$ of $p$ such that  $\sup\limits_U \vert f\vert <\infty.$
\subsection{Statement of the main result}
 Now we are able to state the following
\renewcommand{\thethmspec}{Main Theorem}
\begin{thmspec}
Let  $D=G=E$ and  let $A\subset\partial D, $  $B\subset\partial G$
be  measurable  subsets such that $\mes(A)>0,$
$\mes(B)>0.$
Consider the cross $W:=\X(A,B;D,G).$
Let $M$ be a relatively closed subset of $W$  such that
\begin{itemize}
\item[$\bullet$] $M_a$ is   polar (resp.
discrete)  in $G$ for all $a\in A$ and
 $M^b$
is   polar (resp. discrete)   in $D$  for all $b\in B;$\footnote{ In other words,
$M$ is  polar (resp. discrete) in fibers over $A$ and $B.$}
\item[$\bullet$] $M\cap (A\times B)=\varnothing.$
\end{itemize}
Then there exists a relatively closed pluripolar  subset (resp. an analytic subset)  $\widehat{M}$ of $\widehat{W}$
with the following two properties:
\begin{itemize}
\item[(i)]
The set of end-points of $\widehat{W}\setminus \widehat{M}$ contains  $(W^{\text{o}}\cap W^{\ast})\setminus M.$
\item[(ii)] Let $f:\ W\setminus M \longrightarrow\C$
 be a locally bounded
function  such that
\begin{itemize}
\item[$\bullet$]  for all $a\in A,$
 $f(a,\cdot)|_{G\setminus M_a}$   is
holomorphic and admits the angular limit $f(a,b)$ at all points
$b\in B;$
 \item[$\bullet$] for all $b\in B,$
$f(\cdot,b)|_{D\setminus M^b}$ is holomorphic and admits the
angular limit $f(a,b)$ at all points $a\in A;$
\item[$\bullet$]
$f|_{A\times B}$ is measurable.
\end{itemize}
  Then
there is a unique function
$\hat{f}\in\mathcal{O}(\widehat{W}\setminus \widehat{M})$ such
that $\hat{f}$ admits the angular limit $f$ at  all points of  $(W^{\text{o}}\cap W^{\ast})\setminus M.$
\end{itemize}

Moreover,   if $M=\varnothing,$ then $\widehat{M}=\varnothing.$
\end{thmspec}

\section{Preparatory results}
\subsection{Auxiliary results} First  recall the following
well-known result  (see, for example, \cite{jp1}).
\begin{thm}\label{classical_cross_thm}
Let $D,G$ and $A,\ B$  be open subsets of $\C$ such that $A\subset
D$ and $B\subset G.$ Put $W:=\X(A,B;D,G)$ and
$\widehat{W}:=\widehat{\X}(A,B;D,G).$ Then $W\subset  \widehat{W}$  and every function
$f\in\mathcal{O}_s(W)$ extends uniquely to a function $\hat{f}\in
 \mathcal{O}(\widehat{W}).$
\end{thm}

The following mixed cross theorem has been proved in \cite[Theorem 7.3]{pn2}
(see also  \cite[Theorem 4.2]{nv2} for another proof using the method of holomorphic discs).
\begin{thm}\label{mixed_cross_thm}
 Let $A$ be  a measurable subset of $\partial E$ such that
  $A$ is locally regular.
 Let $G\subset \C$ be an open set and    $B$ an open subset of  $G.$   For $0\leq\delta<1$ put
    $\Omega:=\left\lbrace
   z\in E:\ \omega(z,A,E)<1-\delta\right\rbrace.$
   Let $W:= \X(A,B;\Omega,G)$, $W^{\text{o}}:= \X^{\text{o}}(A,B;\Omega,G),$
      and\footnote{ It will be shown in    Lemma \ref{lem_formula_level_sets} below that
       $\widetilde{\omega}(\cdot,A,\Omega)= \frac{\omega(\cdot,A,E)}{1-\delta}$ on $\Omega,$
where  $\widetilde{\omega}(\cdot,A,\Omega)$ is, in some sense, the ``angular" version of the  harmonic measure.}
         \begin{equation*}
    \widehat{\widetilde{W}}=
     \widehat{\widetilde{\X}}(A,B;\Omega,G):=\left\lbrace (z,w)\in E\times G:\  \frac{\omega(z,A,E)}{1-\delta}
     +\omega
     (w,B,G)<1 \right\rbrace.
     \end{equation*}
Let  $f:\ W\longrightarrow \C$ be   such that
\begin{itemize}
\item[(i)] $f\in\mathcal{O}_s(W^{\text{o}},\C);$
\item[(ii)]  $f$ is  locally  bounded on $W,$  $f|_{A\times B}$ is a 
measurable  function;
\item[(iii)] for all $w\in B,$\footnote{ Since  $A$ is  locally  regular, it follows that $A\subset  \End(\Omega).$}
\begin{equation*}
\lim\limits_{
\mathcal{A}_{\alpha}(a)\ni z\to a }f(z,w)=f(a,w),\qquad a\in A,\
0<\alpha<\frac{\pi}{2}.
\end{equation*}
\end{itemize}
Then there exists a unique function
$\hat{f}\in\mathcal{O}(\widehat{\widetilde{W}})$ such that $\hat{f}=f$ on $\Omega\times B$
and\footnote{    Since  $A$ is  locally  regular, we  have  $A\times G\subset  \End(\widehat{\widetilde{W}}   ).$}
 \begin{equation*}
\lim\limits_{
\mathcal{A}_{\alpha}(a)\ni z\to a  }\hat{f}(z,w)=f(a,w),\qquad a\in A,\
 w\in G,\
0<\alpha<\frac{\pi}{2}.
\end{equation*}
Moreover,   $\vert f\vert_W=\vert \hat{f}\vert_{\widehat{W}}.$
\end{thm}

The next result proved by the authors  in \cite{pn2}
generalizes the work of Gonchar in \cite{go1,go2}.
\begin{thm}\label{Pflug_Nguyen_thm}
We keep the hypotheses and notation of the Main Theorem. Suppose
in addition that $M=\varnothing.$ Then the conclusion of the Main
Theorem  holds for $\widehat{M}=\varnothing.$
\end{thm}

The following two extension theorems are also needed in the
sequel.

\begin{thm}[Chirka \cite{Chi}] \label{Chirka_Thm}
Let $D\subset\C^n$ be a domain and let $\widehat{ D}$ be the
envelope of holomorphy of $D$. Assume that $S$ is a  relatively
closed pluripolar subset of $D$. Then there exists a relatively
closed pluripolar subset $\widehat{ S}$ of $\widehat{ D}$ such
that $\widehat{S}\cap D\subset S$ and $\widehat
{D}\setminus\widehat {S}$ is the envelope of holomorphy of
$D\setminus S$.
\end{thm}

\begin{thm} [Imomkulov--Khujamov \cite{ik}, Imomkulov \cite{im}]\label{Imomkulov_thm}
Let $A$ be  a measurable subset of $\partial E$ with $\mes(A)>0,$ let $M$ be  a relatively closed  subset of $A\times (\C\setminus \overline{E})$
such that  $M_a:=\{w\in\C:\ (a,w)\in M\}$ is  polar (resp. finite) for all $a\in A.$
Then there  exists a  relatively closed pluripolar (resp. analytic) subset $S$ of
$E\times (\C\setminus \overline{E})$  with  the  following  property:

Let  $f:\
(E\cup A)\times E\longrightarrow\C$ be  bounded,  $f|_{E\times E}\in\mathcal{O}(E\times E)$ such
that $\lim\limits_{z\to a,\ z
\in  \mathcal{A}_{\alpha}(a)} f(z,w)=f(a,w)$ for all $ a\in A,
w\in E$ and  $0<\alpha<\frac{\pi}{2}.$ Moreover, assume  that
 the (holomorphic) function $f(a,\cdot)$
extends to a holomorphic function  on $\C\setminus M_a$ for every $a\in A.$ Then  $f|_{E\times E}$
extends holomorphically  to $(E\times \C)\setminus  S.$
\end{thm}
\begin{proof}
This  is a  slightly modified  version of the result  in  \cite{ik}.
In  fact, Imomkulov--Khujamov suppose  that  $f|_{E\times E}$ can be  extended  continuously onto $\overline{E}\times\overline{E}.$
But their proof  still  works under  the hypotheses of  Theorem \ref{Imomkulov_thm}. Consequently,
for each function
$f$  as in  the  statement of the  theorem,
there is   a  relatively closed pluripolar (resp. analytic) subset $S_f$ of
$E\times (\C\setminus \overline{E})$ such that $f$  extends to a  holomorphic  function  on $G_f:= (E\times\C)\setminus S_f$
and that the latter  function does   not extend   holomorphically  across any point of  $ S_f.$
Let $G$ denote the connected component of the interior of 
 $\bigcap_{f}G_f$ that
contains $E^2$ and let $S:=(E\times\C)\setminus G.$  It remains to show
that $S$ is pluripolar (resp. analytic).

Take $(a,b)\in ((A\cap A^\ast)\times\C)\setminus M$. Since $M$ is
relatively closed in $A\times \C$  and  $M_{a}$ is polar, there
exists a smooth curve $\gamma:[0,1]\rightarrow\C\setminus M_{a}$
such that $\gamma(0)=0$, $\gamma(1)=b$. Take an $\epsilon>0$ so
small that
$$
\big (\Delta_a(\epsilon)\times(\gamma([0,1])+\Delta_0(\epsilon))\big )\cap M=\varnothing
$$
and that $V_b:=\Delta_0(\frac{1}{2})\cup(\gamma([0,1])+\Delta_0(\epsilon))$ is  a Jordan  domain.  Consider the cross
\[
Y:=\X(A\cap\Delta_a(\epsilon), \partial V_b\cap
\partial\Delta_0(\frac{1}{2}) ;\Delta_a(\epsilon)\cap E,V_b).
\]
Then $f|_Y$  satisfies the hypotheses  of Theorem  \ref{Pflug_Nguyen_thm}. Consequently,
 we get $\widehat{Y}\subset G_f$ for all $f$ as   in  the  statement of the  theorem.
Hence $\widehat {Y}\subset G$.
Thus $S^{\ast}_{a}\subset M_{a}$ for all $a\in A,$
where  $S^{\ast}_{a}$  is the  non-tangential boundary  layer  of a  pseudoconcave set $S$  (see \cite[p. 358]{im}).
Consequently, by  Lemma  6 and 7  from  \cite{im} (see also Lemma 7 and  8  in \cite{ik}), $S$ is pluripolar (resp. analytic).
\end{proof}
\subsection{Two techniques and their applications}
The technique {\it level sets of (plurisub)harmonic measure} was  introduced by the authors in \cite{pn1}. However,
it turns out that it can be successfully
used
in solving many problems arising from the
theory of separately holomorphic and meromorphic mappings  (see
\cite{pn2,pn3,nv1,nv2}).   For   an open set  $D\subset \C,$  a
subset $A\subset \partial D,$  and $0<\delta<1$   the {\it
$\delta$-level set of the harmonic measure $\omega(\cdot,A,D)$}
is, by definition,
\begin{equation*}
D_{A,\delta}:=\left\lbrace  z\in D:\  \omega(z,A,D)<\delta
\right\rbrace.
\end{equation*}
The technique of level sets  consists in
``replacing" $A$  (resp. $D$) by $D_{A,\delta}$ (resp.
$D_{A,1-\delta}$) for a suitable  $0<\delta<\frac{1}{2}.$

Recall  the following  property of the level sets.
\begin{lem}\label{lem_formula_level_sets}
Let $D$ be  either an empty set or a Jordan domain  such that $E\not\subset D$  and that $D\cup E$
is  a Jordan domain.
For a measurable subset $A$ of $\partial E\cap \partial (D\cup E)$ with $\mes(A)>0$  and $0<\delta <1$ let $\Omega_{\delta}:=E_{A,\delta}\cup D.$
Define the {\rm angular harmonic measure}
\begin{equation*}
\widetilde{\omega}(z,A,\Omega_{\delta}):=
\sup\limits_{u\in\mathcal{U}_{A,\delta}} u(z),  \qquad z\in \Omega_{\delta},
\end{equation*}
where  $\mathcal{U}_{A,\delta}$ is the cone  of all subharmonic functions $u\leq 1$ on $\Omega_{\delta}$ such that
\begin{equation*}
\limsup\limits_{\Omega_{\delta}\cap \mathcal{A}_{\alpha}(\zeta)\ni z\to\zeta}u(z)\leq 0,\qquad
 \zeta\in A\cap A^{\ast},\ 0<\alpha<\frac{\pi}{2}.
\end{equation*}
\begin{itemize}
\item[1)] If $D=\varnothing,$ then  $\widetilde{\omega}(z,A,\Omega_{\delta})=\frac{\omega(z,A,E)}{1-\delta},$  $z\in\Omega_{\delta}.$
\item[2)] If $D$ is a Jordan domain,  then $\widetilde{\omega}(z,A,\Omega_{\delta})\searrow\omega(z,A,E\cup D)$ as  $\delta\searrow 0^{+}.$
\end{itemize}
\end{lem}
\begin{proof}
Part 1) follows  from  \cite[Theorem 4.10]{pn2}.
Part 2)  is a consequence of Part 1).
\end{proof}

The technique {\it conformal mappings} has been introduced by the
second author in \cite{nv2}. This  allows  to reduce the study of
holomorphic extensions on some level sets to the unit disc.

The main idea of the technique of conformal mappings is described  below (see Proposition 5.2 in \cite{nv2} for a proof).
\begin{prop}
\label{prop_conformal}
 Let $A$  be   a measurable subset of $\partial E$  with $\mes(A)>0.$  For $0\leq\delta<1$ put
    $G:=\left\lbrace
   w\in E:\ \omega(w,A,E)<1-\delta\right\rbrace.$
Let $\Omega$ be  an arbitrary connected component of  $G.$
Then
\begin{itemize}
\item[1)] $\End(\Omega) $  is a  measurable subset of $\partial
E$ and $\mes(\End(\Omega))>0.$  Moreover,  $\Omega$ is   a  simply
connected domain.

In virtue of Part 1) and the Riemann mapping theorem,  let $\Phi$
be a conformal mapping of $\Omega$ onto  $E.$

\item[2)] For every  $\zeta\in \End(\Omega),$ there is  $\eta\in\partial
E$ such that
\begin{equation*}
\lim\limits_{ \Omega\cap
\mathcal{A}_{\alpha}(\zeta)\ni z\to  \zeta}\Phi( z)=\eta, \qquad
0<\alpha<\frac{\pi}{2}.
\end{equation*}
$\eta$ is called {\rm the limit of $\Phi$ at the end-point $\zeta$} and it is denoted by $\Phi(\zeta).$
Moreover, $\Phi|_{\End(\Omega)}$  is one-to-one.
\item[3)] Let $f$ be a bounded holomorphic function on  $\Omega,$
$\zeta\in\End(\Omega),$ and $\lambda\in\C$ such that
$\lim\limits_{\Omega\cap
\mathcal{A}_{\alpha}(\zeta)\ni z\to \zeta}f(z)=\lambda$ for some
$0<\alpha<\frac{\pi}{2}.$ Then $f\circ\Phi^{-1}\in\mathcal{O}(E)$
admits the angular limit $\lambda$ at $\Phi(\zeta).$

\item[4)] Let $\Delta$ be a  measurable  subset of  $\End(\Omega)$ such that
$\mes(\Delta)=\mes(\End(\Omega)).$ Put
$\Phi(\Delta):=\{\Phi(\zeta),\ \zeta\in \Delta\},$ where
$\Phi(\zeta)$ is given by Part 2). Then $\Phi(\Delta)$ is a
measurable  subset  of $\partial E$ with  $\mes\big(\Phi(\Delta)  \big)>0$  and
\begin{equation*}
\omega(\Phi(z),\Phi(\Delta),E)=\frac{\omega(z,A,E)}{1-\delta},\qquad
z\in \Omega.
\end{equation*}
\end{itemize}
\end{prop}

As an application of the technique of conformal mappings, we give the following  extended version of Theorem \ref{Imomkulov_thm}.

\begin{thm}\label{new_Imomkulov_thm}
Let $A$ be  a measurable subset of $\partial E$ with $\mes(A)>0.$  For  a given $0\leq\delta<1$ put
    $\Omega:=\left\lbrace
   w\in E:\ \omega(w,A,E)<1-\delta\right\rbrace.$
Let
    $$f:\
\big(\Omega\cup (A\cap \End(\Omega))\big)\times E\longrightarrow\C$$ be a bounded function such
that $f|_{\Omega\times E}$ is holomorphic and $\lim\limits_{z\to a,\ z
\in  \mathcal{A}_{\alpha}(a)} f(z,w)=f(a,w)$ for all $ a\in A\cap \End(\Omega),$
$w\in E$ and  $0<\alpha<\frac{\pi}{2}.$ Suppose in addition that
for every $a\in A\cap \End(\Omega),$ the  function $f(a,\cdot)$
is holomorphic and  it  extends to a holomorphic function  on the whole plane except for a
closed polar (resp. finite)  set of singularities. Then  $f|_{\Omega\times E}$
extends holomorphically  to $(\Omega\times \C)\setminus  S,$ where  $S$
is a  relatively closed pluripolar (resp. analytic) subset of
$\Omega\times \C.$
\end{thm}

 Theorem \ref{Imomkulov_thm}
 is a  special case of the above result for $\delta=0.$
\begin{proof}  We only  treat the  case  where
the set of singularities  of $f(a,\cdot)$ is  closed  polar for $a\in A\cap \End(\Omega).$
Since the  remaining  case  where  these  sets are finite is  analogous, it is  therefore left to the interested  reader.
Using (\ref{eq_elementary_new}) we may   suppose without loss of generality  that  $A$ is locally regular.  Then  $A\subset \End(\Omega).$
Let  $(\Omega_k)_{k\in K}$  be  the family of all  connected components of $\Omega,$  where $K$ is a countable  index set.
By Theorem  4.9  in \cite{pn2},
\begin{eqnarray*}
\End(\Omega)&=&\bigcup\limits_{k\in  K}\End(\Omega_k),\quad  \mes\big(\End(\Omega_k)\cap A\big)= \mes\big(\End(\Omega_k)\big),\\
 & & \End(\Omega_k)\cap\End(\Omega_{k^{'}})=\varnothing\quad
\text{for}\ k\not=k^{'}.
\end{eqnarray*}
By  Proposition \ref{prop_conformal}, we may fix a conformal mapping $\Phi_k$ from $\Omega_k$ onto $E$ for every $k\in K.$
 Put
\begin{equation}\label{eq1_new_Imomkulov_thm}
A_k:=\Phi_k(\End(\Omega_k)\cap A),\  W_k:=(E\cup A_k)\times E,\qquad k\in K.
\end{equation}
Recall from the  hypothesis that for every fixed $w\in E,$ the  holomorphic  function $f(\cdot,w)|_{\Omega}$  is bounded
and that for every $\zeta\in A\cap \End(\Omega),$
\begin{equation*}
\lim\limits_{\Omega\cap \mathcal{A}_{\alpha}(\zeta)\ni z\to \zeta}f(z,w)=f(\zeta,w),\qquad 0<\alpha<\frac{\pi}{2}.
\end{equation*}
Consequently,  Part 3) of   Proposition \ref{prop_conformal}, applied to $f(\cdot,w)|_{\Omega_k}$ with $k\in K,$
implies that for every fixed $w\in E,$ $f(\Phi_k^{-1}(\cdot),w)\in\mathcal{O}(E)$  admits  the angular limit $f(\zeta,w)$  at $\Phi_k(\zeta)$
for all $\zeta\in A\cap \End(\Omega_k).$ By Part 1) of that proposition, we know that $\mes\big( A\cap \End(\Omega_k)\big)>0.$
This  discussion and  the hypothesis
allow us to apply  Theorem \ref{Imomkulov_thm}  to the function
$g_k:\  W_k\longrightarrow\C$  defined by
\begin{equation}\label{eq2_new_Imomkulov_thm}
g_k(z,w):=\begin{cases}
f(\Phi_k^{-1}(z),w),  & (z,w)\in  E\times E\\
f(\Phi_k^{-1}(z),w), & (z,w)\in  A_k\times E
\end{cases},
\end{equation}
where in the second line  we have  used the definition  of $\Phi_k|_{\End(\Omega_k)}$ and  its one-to-one property  proved
by Part 2) of    Proposition \ref{prop_conformal}.
Consequently, we obtain a relatively closed pluripolar set $S_k\subset E\times \C$ such that
$S_k\cap (E\times E) =\varnothing$ and  that
 $g_k|_{E\times E}$
 extends holomorphically to a function $\hat{g}_k\in\mathcal{O}\big((E\times \C)\setminus S_k\big)$
with
\begin{equation}\label{eq3_new_Imomkulov_thm}
 \lim\limits_{
\mathcal{A}_{\alpha}(a)\ni z\to a} \hat{g}_k(z,w)=g_k(a,w),\qquad    (a,w)\in A_k\times E .
\end{equation}
Put
\begin{equation*}
\widehat{\mathcal{W}}_k:=\left\lbrace(\Phi^{-1}_k(z),w),\  (z,w)\in  (E\times\C)\setminus S_k
 \right\rbrace,\qquad  k\in K.
\end{equation*}
Observe that the open sets $ (\widehat{\mathcal{W}}_k)_{k\in K}$  are pairwise disjoint. Moreover,
by (\ref{eq1_new_Imomkulov_thm}),
\begin{multline*}
\bigcup\limits_{k\in K} \widehat{\mathcal{W}}_k= \bigcup\limits_{k\in K}
\left\lbrace (z,w)\in \Omega_k\times \C:\
 (\Phi_k(z),w)\not\in S_k   \right\rbrace\\
= ( \Omega\times \C)\setminus\bigcup\limits_{k\in K}  \left\lbrace (z,w)\in \Omega_k\times \C:\
(\Phi_k(z),w)\in S_k      \right\rbrace=: ( \Omega\times \C)\setminus S.
\end{multline*}
Since $S_k$ is relatively closed pluripolar in $E\times \C$ for $k\in K,$
we see that $S$ is relatively closed pluripolar in $\Omega\times \C.$
Therefore, we  define the desired extension function  $\hat{f}\in\mathcal{O}\big((\Omega\times\C)\setminus S\big)$
  by the formula
\begin{equation*}
\hat{f}(z,w):=\hat{g}_k(\Phi_k(z),w),\qquad   (z,w)\in  \Omega_k,\  k\in K.
\end{equation*}
This, combined with (\ref{eq1_new_Imomkulov_thm})--(\ref{eq3_new_Imomkulov_thm}), implies that  $ \lim\limits_{
 \mathcal{A}_{\alpha}(a)\ni z\to a}\hat{f}(z,w)=f(a,w)$ for all $(a,w)\in A\times E.$
The uniqueness of $\hat{f}$ follows from the one  of $\hat{g}_k,$  $k\in K.$
Hence, the proof of the theorem is complete.
\end{proof}

\subsection{Gluing  theorems}
The following  theorems will be very useful in the next sections when we need to glue different local extensions.

\begin{thm}\label{gluing_thm_1}
Let  $A$ and $\mathcal{N}$ be  measurable subsets of $\partial E$ with $\mes(\mathcal{N})=0.$
 Let $0<\delta<1$   and  $E_{A,\delta}:=\left\lbrace z\in E:\  \omega(z,A,E)<\delta\right\rbrace.$
Suppose that $f\in\mathcal{O}(E_{A,\delta})$  admits the angular limit $0$ at all points
of $(A\cap A^{\ast})\setminus \mathcal{N}.$ Then $f\equiv 0.$
\end{thm}
\begin{proof}
See Theorem 5.4 in \cite{pn2}.
\end{proof}

\begin{thm}\label{minimum_principle}
Let $(\Omega)_{i\in I}$ be a family of open subsets  of an  open set $\Omega\subset \C^n.$
Let  $M$ a relatively closed pluripolar subset (resp. an analytic subset) of $\Omega$ and  $M_i$
 a non pluripolar subset of $\Omega_i$ such that $M\cap M_i=\varnothing,$ $i\in I.$
 Suppose that $f\in\mathcal{O}(\Omega\setminus M)$ and  $f_i\in\mathcal{O}(\Omega_i)$
 satisfy  $f=f_i$ on $M_i$ for $i\in I.$
Then  there exist  a relatively closed pluripolar subset (resp. an analytic subset) $\widehat{M}\subset \Omega$
 and a function
$\hat{f}\in  \mathcal{O}(\Omega\setminus \widehat{M})$  such that
 $\widehat{M}\subset M$  and  $\hat{f}=f$ on $\Omega\setminus M,$
 and that for all $i\in I,$ we have  $\widehat{M}\cap \Omega_i=\varnothing$  and  $\hat{f}=f_i$ on $\Omega_i.$
\end{thm}
\begin{proof} The  case  where  $M$ is a relatively closed pluripolar  subset of $\Omega$ is  not difficult.
The remaining case    where  $M$ is an analytic  subset of $\Omega$ follows  from an easy  application of  Proposition
3.4.5 in \cite{jp1}.
\end{proof}
\begin{thm} \label{gluing_thm_2}
Let $(\Omega_n)_{n=1}^{\infty}$ be an increasing sequence  of open subsets  of an  open set $\Omega\subset \C^n$
such that  $\Omega_n\nearrow\Omega$  as $n\nearrow\infty.$
For every $n\in \N$ let $M_n$ be a
 relatively closed  pluripolar subset  (resp.  an     analytic  subset) of $\Omega_n$
   and  $f_n\in\mathcal{O}(\Omega_n\setminus M_n).$
   Suppose in addition that
 $f_n=f_{n+1} $ on $\Omega_n\setminus (M_n\cup M_{n+1}),$  $n\in\N.$
 Then  there exist  a relatively closed pluripolar subset (resp. an analytic subset)  $M\subset \Omega$
 and a function
$f\in  \mathcal{O}(\Omega\setminus M)$  such that
   $M\cap \Omega_n\subset M_n$  and  $f=f_n$ on $\Omega_n\setminus M_n$ for all $n\in\N.$
\end{thm}
\begin{proof}
It is left  to the interested  reader.
\end{proof}

\section{Extensions through the singularities}
We keep the hypotheses and  notation  of the Main Theorem.
Moreover, we only give the proof for the case where the singular set is {\it fiberwise polar}, that is,
$M_a$ (resp. $M^b$) is polar in $G$  (resp.  $D$) for all $a\in A$  (resp.  $b\in B$).
Since the remaining case  where the singular set is fiberwise discrete is analogous, it is therefore left to the interested reader.

 In this section and  the beginning of  the next one  we assume that

\smallskip

{\bf  $A$  and $B$ are compact sets.}

\smallskip

This  assumption will be  removed  at the end of  the next  section.

 Since $(A\times B)\cap M=\varnothing$ (by the hypothesis), we may find
$N$ points $a_1,\ldots,a_N\in A,$  $N$ numbers $r_1,\ldots,r_N>0,$    $N^{'}$ points
$b_1,\ldots,b_{N^{'}}\in B,$ and $N^{'}$ numbers $s_1,\ldots,s_{N^{'}}>0$ such that
\begin{equation*}
A\subset \bigcup\limits_{k=1}^N\Delta_{a_k}(r_k),\quad B\subset \bigcup\limits_{l=1}^{N^{'}}\Delta_{b_l}(s_l),\quad
M\cap\Big( \bigcup\limits_{k=1}^N\Delta_{a_k}(r_k)\times \bigcup\limits_{l=1}^{N^{'}}\Delta_{b_l}(s_l) \Big)=\varnothing.
\end{equation*}
Put
\begin{equation}\label{subdomains_D_G}
 \widetilde{D}:=D\cap \bigcup\limits_{k=1}^N\Delta_{a_k}(r_k),\quad
\widetilde{G}:=G\cap \bigcup\limits_{l=1}^{N^{'}}\Delta_{b_l}(s_l).
\end{equation}
Then it is  clear that $\X(A,B;\widetilde{D},\widetilde{G})\cap M=\varnothing.$


We introduce  the following  notation.
   For an $a\in A$ (resp. $b\in B$) and $0<r,\delta <1,$ let
 \begin{equation}\label{eq_Dardelta}
 \begin{split}
D_{a,r,\delta}&:=\left\lbrace  z\in D\cap \Delta_a(r):\
\omega(z,A\cap \Delta_a(r),D\cap \Delta_a(r))<\delta
\right\rbrace,\\
 G_{b,r,\delta}&:=\left\lbrace w\in G\cap
\Delta_b(r):\ \omega(w,B\cap \Delta_b(r),G\cap \Delta_b(r))<\delta
\right\rbrace.
\end{split}
\end{equation}

Let  $\Omega$ be an open subset of  $\widehat{W}.$
A point  $(a,b)\in (A\cap A^{\ast})\times G$  (resp.  $(a,b)\in D\times (B\cap B^{\ast})$)  is  said  to be a {\it strong end-point}  of $\Omega$
if there exist $0<r,\delta <1$
 and an open neighborhood $V$ of $b$ (resp. and an open neighborhood  $U$ of $a$) such that
\begin{equation*}
D_{a,r,\delta}\times V\subset \Omega\quad\Big(\text{resp.}\ U\times G_{b,r,\delta}\subset \Omega\ \Big).
\end{equation*}
It is  clear that a strong end-point   of $\Omega$ is  also  an end-point. But the converse statement is  in general  false.

Now, we are in the position to extend $f$ holomorphically through the singular set $M.$
\begin{prop}
\label{prop_extension_2}
For any $a\in A\cap A^{\ast},$   $w\in G,$ there exist $r,\rho,\delta\in ( 0,1) $
   and a
relatively closed pluripolar subset $S\subset D_{a,r,\delta}\times
\Delta_w(\rho)$ with  the  following properties:
\begin{itemize}
\item[1)] $\Delta_w(\rho)\subset G$ and the set $$ T:=\Big( \big (A\cap A^{\ast}\cap
\Delta_a(r)\big) \times \Delta_w(\rho)\Big)\setminus M$$  is  contained  in the  set of strong end-points
of $( D_{a,r,\delta}\times
\Delta_w(\rho))\setminus S.$
\item[2)]  There  is a function $\hat{f}\in\mathcal{O}\Big(
\big(D_{a,r,\delta}\times\Delta_w(\rho)\big )\setminus S\Big) $     which admits
 the angular limit $f$ at all points of $T.$
\end{itemize}
\end{prop}
\begin{proof}
Fix an $a_0\in A\cap A^{\ast}$ and a $w_0\in G$ as in the proposition.
First we
determine  $0<r,\rho,\delta<1$ and then we  will construct  a function $\hat{f}\in\mathcal{O}\Big(
\big(D_{a_0,r,\delta}\times\Delta_{w_0}(\rho)\big )\setminus \widetilde{S}\Big), $
where $ \widetilde{S}$ is a relatively closed pluripolar subset of $D_{a_0,r,\delta}\times\Delta_{w_0}(\rho).$

Since  $M_{a_0}$ is  a relatively closed polar set in $G,$ one may choose
 $\rho>0$  such that $\Delta_{w_0}(\rho)\Subset G$ and
$M_{a_0}\cap\partial\Delta_{w_0}(\rho)=\varnothing$
(cf.~\cite{ArmGar}, Theorem 7.3.9). Take $\rho^-,\rho^+>0$ such that $\rho^-<\rho<\rho^+$,
$\Delta_{w_0}(\rho^+)\Subset G$, and $M_{a_0}\cap\overline
P=\varnothing$, where
\begin{equation*}
P:=\{w\in\C: \rho^-<|w-w_0|<\rho^+\}.
\end{equation*}
Define
\begin{equation*}
\mathcal{G}:=\left\lbrace  w\in \widetilde{G}:\  \omega(w,B,\widetilde{G})<\frac{1}{2}   \right\rbrace.
\end{equation*}
Let $\gamma:[0,1]\to G\setminus M_{a_0}$ be a curve such that
$\gamma(0)\in \mathcal{G}$,
$\gamma(1)\in\partial\Delta_{w_0}(\rho)$. Since $M$ is relatively
closed in $W,$  there exist  $r,t\in (0,1)$ such that
\begin{equation}\label{since_M_is_closed}
 \Delta_{a_0}(r)\cap D\subset \widetilde{D}\quad\text{and}\quad        (A\cap \Delta_{a_0}(r))\times\Big((\gamma([0,1])+\Delta_0(t))\cup
P\Big)\subset W\setminus M.
\end{equation}
Put
\begin{equation*}
V:=\mathcal{G}\cup\big(\gamma([0,1])+\Delta_0(t)\big)\cup
P \end{equation*}
 and consider the cross
\begin{equation*}
Y:=\X(A\cap A^{\ast}\cap \Delta_{a_0}(r),    \mathcal{G};    D_{a_0,r,\frac{1}{2}},V).
\end{equation*}
Using (\ref{subdomains_D_G}) and (\ref{since_M_is_closed})  and the hypotheses on $f$ in the Main Theorem, we
are able to apply Theorem  \ref{Pflug_Nguyen_thm} to the function
$f$ restricted to $\X(A\cap \Delta_{a_0}(r),B;D\cap\Delta_{a_0}(r) ,\widetilde{G}        ).$
Consequently, we obtain  $\widetilde{f}\in \mathcal{O}\big(\widehat{\X}(A \cap\Delta_{a_0}(r),B; D \cap\Delta_{a_0}(r),\widetilde{G} )       \big  )
$ which   admits the angular limit $f$ on
$\X^{\text{o}}(A \cap A^{\ast}\cap \Delta_{a_0}(r),B\cap B^{\ast}; D \cap\Delta_{a_0}(r),\widetilde{G} )       \big  ).
 $
Define
\begin{equation*}
 f_0:=
\begin{cases}
f
  & \text{on}\   (A\cap A^{\ast}\cap \Delta_{a_0}(r))\times V\\
 \widetilde{f} & \text{on}\ D_{a_0,r,\frac{1}{2}}\times
\mathcal{G}
\end{cases}.
\end{equation*}
Then $f_0\in\mathcal{O}_s(Y),$  $f_0|_{ (A\cap A^{\ast}\cap \Delta_{a_0}(r))\times\mathcal{G}}$ is measurable,
and
\begin{multline*}
\lim\limits_{ \mathcal{A}_{\alpha}(\zeta)\ni z\to \zeta}
f_0(z,w)=
f(\zeta,w)=f_0(\zeta,w),\\
 (\zeta,w)\in (A\cap A^{\ast}\cap \Delta_{a_0}(r))\times \mathcal{G},\ 0<\alpha<\frac{\pi}{2}.
\end{multline*}
Consequently, we are able to apply Theorem \ref{mixed_cross_thm} to $f_0$ in order to obtain
 a function $\hat{f}_0$   holomorphic on
 \begin{equation*}
\widehat{Y}=\left\lbrace (z,w)\in   D_{a_0,r,\frac{1}{2}}\times V:\
2\omega\Big(z, A\cap\Delta_{a_0}(r),D\cap \Delta_{a_0}(r)\Big)
 +\omega(w,  \mathcal{G},V)<1    \right\rbrace
\end{equation*}
such that $\hat{f}_0=\widetilde{f}$ on $D_{a_0,r,\delta}\times
\mathcal{G}$ and
\begin{equation}\label{eq_hatf0_prop_extension_2}
\lim\limits_{\mathcal{A}_{\alpha}(\zeta)\ni z\to \zeta}
\hat{f}_0(z,w)=
f(\zeta,w)=:\hat{f}_0(\zeta,w),\ (\zeta,w)\in (A\cap A^{\ast}\cap \Delta_{a_0}(r))\times V,\ 0<\alpha<\frac{\pi}{2}.
\end{equation}
We  have  just  extended  $\hat{f}_0$  to $\widehat{Y}\cup   \big((A\cap A^{\ast}\cap \Delta_{a_0}(r))\times V\big).$
Fix  $s^-,s^+>0$ such that $\rho^-<s^-<\rho< s^+<\rho^+,$ and  consider the annulus
\begin{equation*}
Q:=\{w\in\C: s^{-}<|w-w_0|<s^{+}\}.
\end{equation*}
Let  $\delta$ be such that
\begin{equation*}
0<\delta<\frac{1}{2}\Big(1-\sup\limits_{w\in Q}\omega(w,  \mathcal{G},V)    \Big).
\end{equation*}
Using this and applying Lemma \ref{lem_formula_level_sets}, we see that
 $D_{a_0,r,\delta}\times \overline{Q}\subset \widehat{Y}.$
  Therefore, $\hat{f}_0$ is holomorphic on  $D_{a_0,r,\delta}\times Q$  and continuous on  $D_{a_0,r,\delta}\times \overline{Q}.$ Moreover,
   for any $a\in A\cap A^{\ast}\cap \Delta_{a_0}(r)$ the function
$\hat{f}_0(a,\cdot)$ is holomorphic on $Q$  and  continuous on $\overline{Q}.$ Therefore, by Cauchy formula  we have
\begin{eqnarray*}
 \hat{f}_0(z,w)&=&\frac{1}{2i\pi}\int\limits_{\vert \eta-w_0\vert=s^+}
 \frac{ \hat{f}_0(z,\eta)}{\eta-w}-\frac{1}{2i\pi}\int\limits_{\vert \eta-w_0\vert=s^-}
 \frac{ \hat{f}_0(z,\eta)}{\eta-w}\\
           & =: & \hat{f}^+(z,w)+\hat{f}^-(z,w),\quad
z\in D_{a_0,r,\delta} \cup (A\cap A^{\ast}\cap \Delta_{a_0}(r)) ,\ w\in Q.
\end{eqnarray*}
where $
\hat{f}^+\in\mathcal{O}\Big(  D_{a_0,r,\delta} \times        \Delta_{w_0}(s^+)\Big)$
and
$\hat{f}^-\in\mathcal{O}\Big( D_{a_0,r,\delta}\times (\C\setminus
\overline\Delta_{w_0}(s^-))\Big)$.

Recall from (\ref{eq_hatf0_prop_extension_2}) and the hypotheses that for any $a\in A\cap A^{\ast}\cap \Delta_{a_0}(r)$ the function
$\hat{f}_0(a,\cdot)$ extends holomorphically to $G\setminus M_a.$
Consequently, for any $a\in A\cap A^{\ast}\cap \Delta_{a_0}(r)$ the function
$ \hat{f}^-(a,\cdot)$ extends holomorphically to
$\C\setminus( M_a\cap\overline\Delta_{w_0}(s^-))$. Using (\ref{eq_hatf0_prop_extension_2})
and the above integral formula  for  $\hat{f}^-\in\mathcal{O}\Big( D_{a_0,r,\delta}\times Q\Big),$
we  see that
\begin{equation*}
\lim\limits_{(z,w)\to (\zeta,\eta),\  z\in \mathcal{A}_{\alpha}(\zeta)}
\hat{f}^-(z,w)=\hat{f}^{-}(\zeta,\eta),\qquad (\zeta,\eta)\in (A\cap  A^{\ast}\cap \Delta_{a_0}(r))\times   Q    ,\ 0<\alpha<\frac{\pi}{2}.
\end{equation*}
Now,
we are in the position
to apply
  Theorem \ref{new_Imomkulov_thm} to  $\hat{f}^-.$ Consequently,
there exists a relatively closed pluripolar set
$\widetilde{S}\subset  D_{a_0,r,\delta}\times
\C$  such that $\hat{f}^-$ extends holomorphically to a function
$\overset\approx f{}^-\in\mathcal{O}\big(   (D_{a_0,r,\delta}         \times\C)\setminus \widetilde{S}\big)$.

Since $\hat{f}_0=\hat{f}^++\hat{f}^-$, the function $\hat{f}_0$ extends holomorphically
to a function (still  denoted by)
$\hat{f}_0:=\hat{f}^++\overset\approx f{}^-
\in\mathcal{O}\big(  ( D_{a_0,r,\delta}         \times  \Delta_{w_0}(s^{+}))\setminus \widetilde{S}\big)$.

To prove  Part 1) and Part 2) fix an arbitrary
 $a_1\in A\cap A^{\ast}\cap\Delta_{a_0}(r)$ and $w_1\in \Delta_{w_0}(\rho).$
  Since   $M_{a_1}$
is polar in $G,$ there exists a smooth curve $\alpha:[0,1]\rightarrow\C\setminus M_{a_1}$
such that $\alpha(0)\in \widetilde{G}$ and  $\alpha(1)=w_1.$
Moreover, using  (\ref{subdomains_D_G}) and the hypothesis that   $M$ is a relatively closed subset of $W,$ we may find $r_1>0$ so small that
$\widetilde{V}:= \widetilde{G}\cup(\alpha([0,1])+\Delta_0(r_1))$ is  a Jordan  domain
and that
$$
\Delta_{a_1}(r_1)\Subset \Delta_{a_0}(r),\qquad \big( \Delta_{a_1}(r_1)\times \widetilde{V} \big)\cap  M=\varnothing.
$$
Using this,   
(\ref{subdomains_D_G}), and the hypotheses on $f$ in the Main Theorem, we
are able to apply Theorem  \ref{Pflug_Nguyen_thm} to the function
$f$ restricted to $\X(A\cap \Delta_{a_1}(r_1),B;D\cap\Delta_{a_1}(r_1) ,\widetilde{V}        ).$
Consequently, we obtain  $\widehat{f}_1=\widehat{f}_{(a_1,w_1)}
\in \mathcal{O}\big(\widehat{\X}(A \cap\Delta_{a_1}(r_1),B; D \cap\Delta_{a_1}(r_1),\widetilde{V} )       \big  )
$ which   admits the angular limit $f$ on
$\X^{\text{o}}(A \cap A^{\ast}\cap \Delta_{a_1}(r_1),B\cap B^{\ast}; D \cap\Delta_{a_1}(r_1),\widetilde{V} )       \big  ).
 $
Fix a $w_2\in\mathcal{G}.$ Then  $w_2\in V\cap \widetilde{V}.$
 Choose  $\delta_1,\rho_1>0$ so small  such that
 \begin{equation}\label{eq_choosing_delta_1}
 \delta_1<  1- \omega(w_1, B, \widetilde{V} )
 \end{equation}
 and that
  \begin{equation*}
  D_{a_1,r_1,\delta_1}\times \Delta_{w_2}(\rho_1)\subset   \big( D_{a_0,r,\delta}         \times  \Delta_{w_0}(s^{+}) \big)
         \cap \widehat{\X}(A \cap\Delta_{a_1}(r_1),B; D \cap\Delta_{a_1}(r_1),\widetilde{V} ).
\end{equation*}
 Consequently, using  (\ref{eq_hatf0_prop_extension_2}), we obtain
\begin{multline*}
\lim\limits_{  \mathcal{A}_{\alpha}(\zeta)\ni z\to \zeta}
\hat{f}_0(z,w)=\lim\limits_{  z\to \zeta,\  z\in \mathcal{A}_{\alpha}(\zeta)}
\hat{f}_1(z,w)=f(\zeta,w),\\
                \zeta\in A\cap A^{\ast}\cap \Delta_{a_1}(r_1),\ w\in \Delta_{w_2}(\rho_1),    0<\alpha<\frac{\pi}{2}.
\end{multline*}
By Theorem \ref{gluing_thm_1}, $\hat{f}_0=\hat{f}_1$ on $ D_{a_1,r_1,\delta_1}\times \Delta_{w_2}(\rho_1).$
By shrinking $\rho_1$ (if necessary) and by using the fact that $w_1,w_2\in \widetilde{V}$ and estimate (\ref{eq_choosing_delta_1}),
 we deduce from the latter
identity that
\begin{equation}\label{last_eq_prop_extension_2}
\hat{f}_0=\hat{f}_1\quad\text{on}\ \big(D_{a_1,r_1,\delta_1}\times \Delta_{w_1}(\rho_1)\big)\setminus  \widetilde{S}.
\end{equation}
Now   we  are in the position to apply
 Theorem \ref{minimum_principle} to  $\hat{f}_0
\in\mathcal{O}\big(  ( D_{a_0,r,\delta}         \times  \Delta_{w_0}(s^{+}))\setminus \widetilde{S}\big)$   and to
 the  family of   functions $\big(\widehat{f}_{(a_1,w_1)}\big)$ with
  $a_1\in A\cap A^{\ast}\cap\Delta_{a_0}(r)$ and $w_1\in \Delta_{w_0}(\rho)\setminus M_{a_1}.$
 Consequently, we obtain  the  relatively  closed pluripolar     subset $S\subset D_{a_0,r,\delta}\times\Delta_{w_0}(\rho)
 $ satisfying Part 1) and  the   function $\hat{f}\in\mathcal{O}\Big(
\big(D_{a_0,r,\delta}\times\Delta_{w_0}(\rho)\big )\setminus S\Big). $
 Part 2)  follows from (\ref{last_eq_prop_extension_2}).
 \end{proof}
 The  role  of  strong end-points is  illustrated  by the  following uniqueness theorem.

 \begin{thm}\label{uniqueness}
 Let $f\in\mathcal{O}(\Omega),$  where $\Omega$ is   a  subdomain  of $\widehat{W}.$
 Suppose that there exist $a_0\in A\cap A^{\ast},$  $r>0$ and  an open subset
 $V\subset G$  such that  $(A\cap A^{\ast}\cap \Delta_{a_0}(r))\times V$ are contained in the set of strong end-points of $\Omega$
 and that  the angular limit of $f$ at all points of  $(A\cap A^{\ast}\cap \Delta_{a_0}(r))\times V$ equals $0.$
 Then $f\equiv 0.$
 \end{thm}
 \begin{proof}
 Applying Theorem \ref{gluing_thm_1}  to $f$ restricted  to an open set of the form $D_{a_0,r_0,\delta}\times U\subset \Omega$
 for suitable $r_0,\delta>0$ and $U\subset G,$ the   theorem follows.
 \end{proof}
\section{Proof of the Main Theorem} 
We keep the notation in the previous section. Moreover,
we introduce some  new notation.
 For any $\zeta\in A,$    $r,  R\in (0,1)$  with  $\Delta_0(R)\cap \widetilde{G}\not=\varnothing,$ let
 \begin{equation*}
 W_{\zeta,r, R}:= \X\big(A\cap \Delta_{\zeta}(r),B;  D\cap \Delta_{\zeta}(r),
\Delta_0(R)\cup\widetilde{G}\big).
\end{equation*}
Similarly,
 for any $\eta\in B,$   $r,  R\in (0,1)  $  with  $\Delta_0(R)\cap \widetilde{D}\not=\varnothing,$ put
 \begin{equation*}
 W_{\eta,r, R}:= \X\big(A,B\cap \Delta_{\eta}(r);
\Delta_0(R)\cup\widetilde{D},  G\cap \Delta_{\eta}(r) \big).
\end{equation*}
For  any $R\in (0,1)$  with  $\Delta_0(R)\cap \widetilde{D}\not=\varnothing$
and  $\Delta_0(R)\cap \widetilde{G}\not=\varnothing,$ put
 \begin{equation*}
 W_{ R}:= \X\big(A,B;
\Delta_0(R)\cup\widetilde{D}, \Delta_0(R^{'})\cup\widetilde{G}  \big).
\end{equation*}
Fix a sequence $(\delta_n:=\frac{1}{2^n})_{n=1}^{\infty}.$

The proof is divided into several steps. In the  first three steps  $A$ and $B$ are supposed  to be compact.

\smallskip

\noindent {\bf Step 1.} { \it For any $\zeta\in A\cap A^{\ast}$  and $R\in (0,1)$  with  $\Delta_0(R)\cap \widetilde{G}\not=\varnothing,$
 there exists $r\in (0,1)$ and   a relatively closed pluripolar subset $\widehat{S}$ of
$
\widehat{W}_{\zeta,r,R}$
and
a function $\hat{f}\in \mathcal{O}(\widehat{W}_{\zeta,r,R} \setminus \widehat{S}) $  with  the following  properties:
\begin{itemize}
\item[$\bullet$] $ (W^{\text{o}}_{\zeta,r,R}\cap W^{\ast}_{\zeta,r,R})\setminus M $ is  contained in the set
of strong  end-points  of $\widehat{W}_{\zeta,r,R} \setminus \widehat{S}.$
\item[$\bullet$] $\hat{f}$  admits the  angular limit
$f$ at all points  of
$ (W^{\text{o}}_{\zeta,r,R}\cap W^{\ast}_{\zeta,r,R} )\setminus M.$
\end{itemize}
}

Applying Proposition
\ref{prop_extension_2}
to the points $\zeta\in A\cap A^{\ast}$ and $ w\in \overline{\Delta}_0(R)$ and using the compactness of $ \overline{\Delta}_0(R),$
 we  find  $r,\delta\in (0,1),$ $p\in \N,$   and for  any $j\in\{1,\ldots,p\},$ a point  $w_j\in  \overline{\Delta}_0(R),$  a number
  $\rho_j>0,$  a  relatively  closed  subset $S_j\subset D_{\zeta,r,\delta}\times  \Delta_{w_j}(\rho_j),$
and  a function  $\hat{f}_j\in\mathcal{O}\big(   (D_{\zeta,r,\delta}\times  \Delta_{w_j}(\rho_j)) \setminus S_j \big)$ such that
 \begin{itemize}
 \item[$\bullet$] $\overline{\Delta}_0(R)\subset\bigcup\limits_{k=1}^p \Delta_{w_k}(\rho_k);$
 \item[$\bullet$]
$\hat{f}_j$  admits  the angular limit $f$ at all points  of
$\Big (\big ( A\cap A^{\ast}\cap
\Delta_{\zeta}(r)  \big)\times \Delta_{w_j}(\rho_j)\Big)\setminus M.$
\end{itemize}
Using this  we are able  to apply  Theorem \ref{uniqueness}. Consequently,
 $$\hat{f}_i=\hat{f}_j\quad\text{on}\
 \big( D_{\zeta,r,\delta}\times  (\Delta_{w_i}(\rho_i)\cap \Delta_{w_j}(\rho_j)\cap \Delta_0(R) ) \big) \setminus (S_i\cup S_j).$$
 Therefore, we obtain an  $\tilde{f}\in \mathcal{O}\big(  ( D_{\zeta,r,\delta}\times\Delta_0(R))\setminus S^{'}\big),$ where
 $\tilde{f} =\hat{f}_j$ on $ (D_{\zeta,r,\delta}\times  \Delta_{w_j}(\rho_j)) \setminus S^{'}$
 and $S^{'}:=\bigcup\limits_{j=1}^p S_j$  is relatively  closed  pluripolar set.  Moreover, $\Big(\big ( A\cap A^{\ast}\cap
\Delta_{\zeta}(r)  \big)\times \Delta_0(R)\Big)\setminus M$  is contained  in the set of strong end-points of
$ \big ( D_{\zeta,r,\delta}\times\Delta_0(R)\big)\setminus S^{'}$  and $\tilde{f}$  admits  the angular limit $f$ at all points  of
the former set.

On the other hand,   applying Theorem  \ref{Pflug_Nguyen_thm} to the function
$f$ restricted to $\X(A\cap \Delta_{\zeta}(r),B;D\cap\Delta_{\zeta}(r) ,\widetilde{G}        ),$
 we obtain  $\overset\approx f{}
\in \mathcal{O}\big(\widehat{\X}(A \cap\Delta_{\zeta}(r),B; D \cap\Delta_{\zeta}(r),\widetilde{G} )       \big  )
$ which   admits the angular limit $f$ on
$\X^{\text{o}}(A \cap A^{\ast}\cap \Delta_{\zeta}(r),B\cap B^{\ast}; D \cap\Delta_{\zeta}(r),\widetilde{G} )       \big  ).
 $

Next, we fix an $n_0$  such that $\delta_{n_0}<\delta.$ For $s\in (0,1)$ let
$\widetilde{G}_s:=\{w\in \widetilde{G}:\  \omega(w,B,\widetilde{G})<s\}.$
  For all $n\geq n_0$ let
  \begin{equation*}
  W_n:=\X\Big( D_{\zeta,r,\delta_n}, \widetilde{G}_{\delta_n};  D_{\zeta,r,1-\delta_n},
    \Delta_0(R)\cup  \widetilde{G}_{1-\delta_n}\Big).
    \end{equation*}
Define $f_n:\  W_n\setminus S^{'}\rightarrow\C$  as  follows
\begin{equation}\label{eq_fn}
f_n:=
\begin{cases}
 \tilde{f}, &\text{on}\ \Big(D_{\zeta,r,\delta_n}\times \big ( \Delta_0(R)\cup  \widetilde{G}_{1-\delta_n}\big)\Big)\setminus S^{'}\\
 \overset\approx f{}, &\text{on}\  D_{\zeta,r,1-\delta_n}\times  \widetilde{G}_{\delta_n}
\end{cases};
\end{equation}
  here  we  have  applied Theorem \ref{gluing_thm_1} in order to show that $ \tilde{f}=\overset\approx f{}$
  on the overlapping  set. Clearly, $f_n\in\mathcal{O}( W_n\setminus S^{'}).$
  Therefore, applying Theorem \ref{classical_cross_thm}   and Theorem  \ref{Chirka_Thm} to  $ W_n\setminus S^{'},$   we obtain
 a  relatively closed
pluripolar subset  $ \widehat{S}_n$ of  $ \widehat{W}_n$ with  $ \widehat{S}_n\cap W_n\subset S^{'}$ and   a function
$\hat{f}_n\in \mathcal{O}( \widehat{W}_n\setminus \widehat{S}_n)$
 with  $\hat{f}_n=f_n$ on $  W_n\setminus S^{'}.$  Now, using  Lemma \ref{lem_formula_level_sets}, we define
 \begin{eqnarray*}
  X_n&:=&\X\big(  A\cap A^{\ast}\cap\Delta_{\zeta}(r), B\cap B^{\ast}  ;  D_{\zeta,r,1-\delta_n},
    \Delta_0(R)\cup  \widetilde{G}_{1-\delta_n}\big),\\
    \widehat{X}_n&:=&\left\lbrace (z,w)\in  D_{\zeta,r,1-\delta_n} \times ( \Delta_0(R)\cup  \widetilde{G}_{1-\delta_n}):\
    \widetilde{\omega}\big(z, A\cap A^{\ast}\cap\Delta_{\zeta}(r),  D_{\zeta,r,1-\delta_n} \big ) \right.\\
   &\qquad& \left. +   \widetilde{\omega}\big(w, B\cap B^{\ast},   \Delta_0(R)\cup  \widetilde{G}_{1-\delta_n}\big )<1 \right\rbrace.
\end{eqnarray*}
 Then it follows from (\ref{eq_fn}) that  $\hat{f}_n$  restricted
 to $\widehat{X}_n\setminus \widehat{S}_n,$
     admits the angular limit $f$  at all points
 of $X^{\text{o}}_n$  and  the latter set is  contained  in the set of strong end-points
    of $\widehat{X}_n\setminus \widehat{S}_n.$
    Therefore,  applying  Theorem \ref{uniqueness}    we  see that
    $\hat{f}_n=\hat{f}_{n+1}$ on $\widehat{X}_n\setminus (\widehat{S}_n\cup\widehat{S}_{n+1}). $
    Moreover, using   Theorem  \ref{minimum_principle} we may assume that
$\widehat{S}_{n+1}\cap\widehat{X}_n \subset \widehat{S}_n .$
    Next, we will show  that $\widehat{X}_n\nearrow \widehat{W}_{\zeta,r,R}$ as $n\nearrow\infty.$ To see this
   it suffices to observe by Lemma \ref{lem_formula_level_sets} that
    \begin{eqnarray*}
    \omega\big (\cdot, A\cap A^{\ast}\cap\Delta_{\zeta}(r), D_{\zeta,r,1-\delta_n}\big)&\searrow&
   \omega(\cdot, A\cap A^{\ast}, D\cap\Delta_{\zeta}(r)),\\
   \omega\big (\cdot , B\cap B^{\ast}, \Delta_0(R)\cup  \widetilde{G}_{1-\delta_n}\big)&\searrow&
     \omega(\cdot, B\cap B^{\ast},\widetilde{G}),
     \end{eqnarray*}
      when $n\nearrow\infty.$
      Now  we are in the position to apply Theorem    \ref{gluing_thm_2} to the functions
   $\hat{f}_{n}\in\mathcal{O}(\widehat{X}_n\setminus \widehat{S}_n)$ for $n\geq n_0.$
  Consequently, we  obtain the  desired relatively  closed  pluripolar subset $\widehat{S}$ of $ \widehat{W}_{\zeta,r,R} $
  and the  desired  extension  function $\hat{f}.$
      This finishes Step 1.

  \smallskip

\noindent {\bf Step 2.} { \it
For any $R\in (0,1)$  such that  $\Delta_0(R)\cup \widetilde{D}$ and
 $\Delta_0(R)\cup \widetilde{G}$ are Jordan domains,
  there exist a relatively closed pluripolar subset $\widehat{S}$ of
$
\widehat{W}_{R}
$
and
a function $\hat{f}\in \mathcal{O}(\widehat{W}_{R} \setminus \widehat{S}) $
such that  the set $( W^{\text{o}}_{R} \cap W^{\ast}_{R})\setminus M $  is  contained  in the  set
of strong end-points of $\widehat{W}_{R} \setminus \widehat{S}$  and that $\hat{f}$  admits the angular limit
$f$ at all points  of the  former set.
}

\smallskip

Choose a sequence of closed subsets $(\widetilde{A}_m)_{m=1}^{\infty}$ (resp. $(\widetilde{B}_m)_{m=1}^{\infty}$)  of $\partial D$
(resp.  $\partial G$)  such that
\begin{equation}\label{eq1_Step2}
\begin{split}
\mes(\widetilde{A}_m)&>0, \ \widetilde{A}_m\subset \widetilde{A}_{m+1}\subset A\cap A^{\ast},\  \mes\Big( A\setminus
 \bigcup_{m=1}^{\infty}\widetilde{A}_m\Big)=0,\\
  \mes(\widetilde{B}_m)&>0,\ \widetilde{B}_m\subset \widetilde{B}_{m+1}\subset  B\cap B^{\ast} ,\ \mes\Big( B\setminus \bigcup_{m=1}^{\infty}
  \widetilde{B}_m\Big)=0.
  \end{split}
\end{equation}
Let
\begin{equation}\label{eq2_Step2}
\mathcal{W}_m:= \X(\widetilde{A}_m,\widetilde{B}_m;\Delta_0(R)\cup  \widetilde{D},\Delta_0(R)\cup  \widetilde{G}),
\quad \widehat{\mathcal{W}}_m := \widehat{\X}
(\widetilde{A}_m,\widetilde{B}_m;\Delta_0(R)\cup  \widetilde{D},\Delta_0(R)\cup  \widetilde{G}).
\end{equation}
First, we will show that for every $m$  there exist a relatively closed pluripolar subset $\widehat{\mathcal{S}}_m$ of
$
\widehat{\mathcal{W}}_{m}
$
and
a function $\widetilde{f}_m\in \mathcal{O}(\widehat{\mathcal{W}}_{m} \setminus \widehat{\mathcal{S}}) $
such that  the set $( \mathcal{W}^{\text{o}}_{m} \cap \mathcal{W}^{\ast}_{m})\setminus M $  is  contained  in the  set
of strong end-points of $\widehat{\mathcal{W}}_{m} \setminus \widehat{\mathcal{S}}_m$  and that $\widetilde{f}_m$  admits the angular limit
$f$ at all points  of the  former set.
For this  purpose fix  an $m\in \N.$

Applying  Step 1 and  using  a compactness argument with respect to $\widetilde{A}_m$ we may find  $K$ points $\zeta_1,\ldots,\zeta_K\in A\cap A^{\ast}$ and
$K$ numbers $r_1,\ldots,r_K>0$ with the following properties:
\begin{itemize}
\item[$\bullet$]
    $\widetilde{A}_m\subset \bigcup\limits_{k=1}^K\Delta_{\zeta_k}(r_k)$ and  $ D\cap \bigcup\limits_{k=1}^K\Delta_{\zeta_k}(r_k) \subset \widetilde{D};$
\item[$\bullet$] for every  $1\leq k\leq K,$ there are
   a relatively closed pluripolar subset $S_k$ of $\widehat{W}_{\zeta_k,r_k,R}$ and
a function $\hat{g}_k\in \mathcal{O}\big(\widehat{W}_{\zeta_k,r_k,R} \setminus S_k\big)$ such that
 the set $ (W^{\text{o}}_{\zeta_k,r_k,R} \cap W^{\ast}_{\zeta_k,r_k,R})\setminus M$ is  contained
 in the  set of strong end-points of  $\widehat{W}_{\zeta_k,r_k,R}\setminus S_k$
and that $\hat{g}_k$   admits the angular limit
$f$ at all points  of the former set.
\end{itemize}
Similarly, using  Step 1 again but exchanging the role between $A$ and $B$ (resp. $D$ and $G$),
we may find  $L$ points $\eta_1,\ldots,\eta_{L}\in B\cap B^{\ast}$ and
$L$ numbers $s_1,\ldots,s_L>0$  
   with the following  properties:
\begin{itemize}
\item[$\bullet$]
    $\widetilde{B}_m\subset \bigcup\limits_{l=1}^{L}\Delta_{\eta_{l}}(s_{l})$ and
    $ G\cap \bigcup\limits_{l=1}^{L}\Delta_{\eta_{l}}(s_{l}) \subset \widetilde{G};$
    \item[$\bullet$] for every  $1\leq l\leq L,$ there are
      a relatively closed pluripolar subset $T_{l}$ of $\widehat{W}_{\eta_{l},s_{l},R} $
 and a function  $\hat{h}_{l}\in \mathcal{O}\big(\widehat{W}_{\eta_{l}, s_{l},   R}\setminus T_{l}\Big)$
such that the set $ (W^{\text{o}}_{\eta_{l},s_{l},R} \cap W^{\ast}_{\eta_{l},s_{l},R})\setminus M$ is  contained
 in the  set of strong end-points of  $\widehat{W}_{\eta_{l},s_{l},R}\setminus T_{l}$
and that $\hat{h}_{l}$
 admits the angular limit
$f$ at all points  of the former set.
\end{itemize}
Put $S:=\bigcup\limits_{k=1}^K S_k$ and $ T:= \bigcup\limits_{l=1}^{L} T_{l}.$
For every $n\geq 1$  let
\begin{eqnarray*}
A_n  &:=& \bigcup\limits_{k=1}^K D_{\zeta_k,r_k,\delta_n},\quad \quad B_{n}:=
\bigcup\limits_{l=1}^{L} G_{\eta_{l},s_{l},\delta_n},\\
D_n&:=&\left\lbrace       z\in \Delta_0(R)\cup \widetilde{D}:\ \omega\Big( z, \widetilde{A}_m,  \Delta_0(R)\cup \widetilde{D}      \Big)
<1-\delta_n             \right\rbrace,\\
G_{n}  &:=&\left\lbrace       w\in \Delta_0(R)\cup \widetilde{G}:\ \omega\Big( w, \widetilde{B}_m,  \Delta_0(R)\cup \widetilde{G}      \Big)
<1-\delta_n             \right\rbrace,\\
W_{n}  &:=& \X\Big(A_{n},
B_{n};D_{n},G_{n}\Big),\qquad X_n:=\X(\widetilde{A}_m\cap \widetilde{A}_m^{\ast},\widetilde{B}_m\cap \widetilde{B}_m^{\ast};D_n,G_n),\\
\widehat{\widetilde{X}}_n &:=& \left\lbrace (z,w)\in D_n\times G_n:\  \widetilde{\omega}(z,\widetilde{A}_m\cap \widetilde{A}_m^{\ast},D_n)+
\widetilde{\omega}(w,\widetilde{B}_m\cap \widetilde{B}_m^{\ast},G_n)<1
         \right\rbrace    ,
\end{eqnarray*}
where  in the last line we  can apply  Lemma \ref{lem_formula_level_sets}  since $\Delta_0(R)\cup \widetilde{D}$
and $  \Delta_0(R)\cup \widetilde{G}$  are Jordan  domains.
  Applying Theorem \ref{uniqueness} and Theorem \ref{gluing_thm_2}, we may glue $(\hat{g}_k)_{k=1}^K$ together
 in order to  define the function
$g_{n}:\ (A_n\times G_n)\setminus S\longrightarrow\C$  as follows
\begin{equation*}
g_{n} :=
\hat{g}_k
  \qquad \text{on}\  (D_{\zeta_k,r_k,\delta_n}\times  G_{n})  \setminus S.
\end{equation*}
Similarly,  we may glue $(\hat{h}_{l})_{l=1}^{L}$ together
 in order to  define the function
$h_{n}:\ (D_n\times B_n)\setminus T\longrightarrow\C$  as follows
\begin{equation*}
h_{n} :=
\hat{h}_{l}
  \qquad \text{on}\ (D_n\times  G_{\eta_{l},s_{l},\delta_n} ) \setminus T.
\end{equation*}
Finally, we  glue $g_n$ and $h_n$  together in order to  define the function
$f_{n}:\ W_n\setminus (T\cup S)\longrightarrow\C$  as follows
\begin{equation}\label{eq_new_fn}
f_{n} :=
\begin{cases}
g_n
  &\text{on}\ (A_n\times  G_{n})  \setminus S\\
 h_n, & \text{on}\ (D_{n}\times  B_{n}) \setminus T
\end{cases}.
\end{equation}
The  remaining  part of the proof follows  along the same lines as in Step 1.
Applying Theorem \ref{classical_cross_thm}   and Theorem  \ref{Chirka_Thm} to  $ W_n\setminus( S\cup T),$   we obtain
 a  relatively closed
pluripolar subset  $ \widehat{S}_n$ of  $ \widehat{W}_n$ with  $ \widehat{S}_n\cap W_n\subset (S\cup T)$ and   a function
$\hat{f}_n\in \mathcal{O}( \widehat{W}_n\setminus \widehat{S}_n)$
 with  $\hat{f}_n=f_n$ on $  W_n\setminus (S\cup T).$
 In particular, it follows from (\ref{eq_new_fn}) that  $\hat{f}_n$  restricted
 to $\widehat{\widetilde{X}}_n\setminus \widehat{S}_n,$
     admits the angular limit $f$  at all points
 of $X^{\text{o}}_n\setminus M.$ Here observe that    the latter set is  contained  in the set of strong end-points
    of $\widehat{\widetilde{X}}_n\setminus \widehat{S}_n.$
    Therefore,  applying  Theorem \ref{uniqueness}    we  see that
    $\hat{f}_n=\hat{f}_{n+1}$ on $\widehat{\widetilde{X}}_n\setminus (\widehat{S}_n\cup\widehat{S}_{n+1}). $
    Moreover, using   Theorem  \ref{minimum_principle} we may assume that
$\widehat{S}_{n+1}\cap\widehat{\widetilde{X}}_n \subset \widehat{S}_n .$
    Next, an application of Lemma \ref{lem_formula_level_sets} gives that
     $\widehat{\widetilde{X}}_n\nearrow \widehat{\mathcal{W}}_{m}$ as $n\nearrow\infty.$
      Now  we are in the position to apply Theorem    \ref{gluing_thm_2} to the functions
   $\hat{f}_{n}\in\mathcal{O}(\widehat{\widetilde{X}}_n\setminus \widehat{S}_n)$ for $n\geq 2.$
  Consequently, we  obtain the  desired relatively  closed  pluripolar subset $\widehat{\mathcal{S}}_m$ of $ \widehat{\mathcal{W}}_{m} $
  and the  desired  extension  function $\widetilde{f}_m.$

 Using   (\ref{eq1_Step2})--(\ref{eq2_Step2})  and     (\ref{eq_elementary}) we see that
$\widehat{\mathcal{W}}_m\nearrow \widehat{W}_R$ as $m\nearrow\infty.$
Therefore,   applying
 Theorem \ref{uniqueness},  Theorem  \ref{minimum_principle} and  Theorem    \ref{gluing_thm_2},    Step 2  follows.

  \smallskip

\noindent {\bf Step 3.} { \it Completion of the case where $A$ and $B$ are compact.}

\smallskip

Fix  $R\in (0,1)$   such that $\Delta_0(R)\cup \widetilde{D}\not=\varnothing$
and  $\Delta_0(R)\cup \widetilde{G}\not=\varnothing$  are Jordan domains. Choose a sequence
$(R_n)_{n=1}^{\infty}$   such that
$R_n >R,$  $R_n\nearrow 1$ as $n\nearrow\infty.$
For $n\geq  1$ put $
W_{n}  :=  \X\Big(A,B; \Delta_0(R_n)\cup \widetilde{D}    ,\Delta_0(R_n)\cup \widetilde{G}\Big).
$
Applying the result of Step  2,  we may find, for every $n\geq 1,$   a relatively closed pluripolar subset $\widehat{M}_n$ of $\widehat{W}_{n}$
and
a function $\hat{f}_n\in \mathcal{O}(\widehat{W}_{n} \setminus \widehat{M}_n) $ such that
 $(W_{n}^{\ast}\cap W^{\text{o}}_n)\setminus M$  is contained  in the set of strong end-points of $\widehat{W}_{n} \setminus \widehat{M}_n$ and that
$\hat{f}_n$  admits the angular limit
$f$ at all points  of the  former set.
Since  $\widehat{W}_{n}\nearrow \widehat{W}$ as $n\nearrow \infty,$ we  conclude  this  step by  applying
 Theorem \ref{uniqueness},  Theorem  \ref{minimum_principle} and  Theorem    \ref{gluing_thm_2} as in Step 2.


\smallskip

\noindent {\bf Step 4.} { \it The general case.}

\smallskip

Choose a sequence of closed subsets $(A_n)_{n=1}^{\infty}$ (resp. $(B_n)_{n=1}^{\infty}$)  of $\partial D$
(resp.  $\partial G$)  such that
\begin{eqnarray*}
\mes(A_n)&>&0, \ A_n\subset A_{n+1}\subset A,\  \mes\Big( A\setminus \bigcup_{n=1}^{\infty}A_n\Big)=0,\\
  \mes(B_n)&>&0,\ B_n\subset B_{n+1}\subset B,\ \mes\Big( B\setminus \bigcup_{n=1}^{\infty}B_n\Big)=0.
\end{eqnarray*}
Let $W_n:= \X(A_n,B_n;D,G).$
Applying  the hypotheses  to $f|_{W_n\setminus M}$ for $n\geq 1,$
we obtain a relatively closed pluripolar subset $\widehat{M}_n$ of $\widehat{W}_n$ and
a function $\hat{f}_n\in\mathcal{O}(\widehat{W}_n\setminus \widehat{M}_n)$
such that  $(W_{n}^{\ast}\cap W^{\text{o}}_n)\setminus M$  is contained  in the set of strong end-points of $\widehat{W}_{n} \setminus \widehat{M}_n$ and that
$\hat{f}_n$  admits the angular limit
$f$ at all points  of the  former set.   Using (\ref{eq_elementary}) we see that
$\widehat{W}_n\nearrow \widehat{W}$ as $n\nearrow\infty.$
Therefore, arguing as  at the end of the previous step,    Step 4  follows.
 \hfill  $\square$

%
%
%
%

%
%
%
%
%

\end{document}